\documentclass{article}
\usepackage{arxiv}
\usepackage{lineno,hyperref}
\modulolinenumbers[1]
\usepackage{amsmath}
\usepackage{graphicx}
\usepackage{amsfonts}
\usepackage{float}
\usepackage{verbatim}
\usepackage{amsthm}
\usepackage{amssymb}
\usepackage{latexsym}
\usepackage{epstopdf}
\usepackage{color}
\usepackage{a4}
\usepackage{amsmath}
\allowdisplaybreaks[4]

\makeatletter
\@addtoreset{equation}{section}

\newcommand{\EE}{\mathbb{E}}
\newcommand{\PP}{\mathbb{P}}
\newcommand{\RR}{\mathbb{R}}

\newcommand{\NN}{\mathbb{N}}

\newcommand{\ind}{\mathbf{1}}

\def\F{{\cal F}}
\def\G{{\cal G}}

\newtheorem{assp}{Assumption}
\newtheorem{expl}{Example}
\newtheorem{coro}{Corollary}
\newtheorem{rmk}{Remark}
\newtheorem{theorem}{Theorem}
\newtheorem{lemma}{Lemma}

\newcommand{\eproof}{\indent\vrule height6pt width4pt depth1pt\hfil\par\medbreak}

\def\a{\alpha}
   
 \def\k{\kappa} \def\l{\lambda} \def\m{\mu} 
  \def\r{\rho}

\def\D{\Delta}

\def\1{\oslash} \def\2{\oplus} \def\3{\otimes} \def\4{\ominus}
\def\5{\circ} \def\6{\odot} \def\7{\backslash} \def\8{\infty}
\def\9{\bigcap} \def\0{\bigcup} \def\+{\pm} \def\-{\mp}

\title{Truncated Euler-Maruyama method for time-changed stochastic differential equations with super-linear state variables and H\"older's continuous time variables}

\author{Xiaotong Li ~~ Wei Liu ~~ Tianjiao Tang \\ Department of Mathematics\\ Shanghai Normal University, Shanghai, 200234, China}

\begin{document}
\maketitle

\begin{abstract}
An explicit numerical method is developed for a class of non-autonomous time-changed stochastic differential equations, whose coefficients obey H\"older's continuity in terms of the time variables and are allowed to grow super-linearly in terms of the state variables. The strong convergence of the method in the finite time interval is proved and the convergence rate is obtained. Numerical simulations are provided.
\medskip  \par \noindent
{\small\bf Key words}: explicit numerical method, highly non-linear coefficients, time-changed processes, stochastic differential equations, strong convergence.
\medskip  \par \noindent
{\small\bf 2010 MSC}: 65C30, 60H10
\end{abstract}

\section{Introduction}
\par \noindent
Time-changed processes and time-changed stochastic differential equations (SDEs), as important mathematical tools to describe subdiffusions, have been broadly investigated in recent decades \cite{Chen2017,Kobayashi2011,Liu2021,Magdziarz2009,MS2004,NN2018,UHK2018}. As the explicit forms of the underlying solutions to time-changed SDEs are rarely obtained, numerical methods become extremely important when time-changed SDEs are applied in practice.
\par
For time-changed SDEs of the following form
\begin{equation}
\label{eq:duatcSDE}
dY(t) = f(E(t),Y(t))dE(t) + g(E(t),Y(t))dB(E(t)),
\end{equation}
the Euler-Maruyama (EM) method was proposed in \cite{JK2016}, in which both the coefficients $f$ and $g$  were required to satisfy the global Lipschitz condition in terms of the state variable. To our best knowledge, \cite{JK2016} is the paper that investigates the finite time convergence of the numerical method for time-changed SDEs. When the constraint on $f$ is released to the one-sided Lipschitz condition for \eqref{eq:duatcSDE}, the semi-implicit EM method was studied in \cite{DL2020online}. Furthermore, the truncated EM method was discussed in \cite{LMTW2020} when some super-linear terms are allowed to appear in both the coefficients $f$ and $g$.  All the three works mentioned above employed the duality principle that was developed in \cite{Kobayashi2011}. Briefly speaking, the duality principle states that the solution to \eqref{eq:duatcSDE} can be represented by $X(E(t))$, where $X(t)$ is the solution to the classical SDE
\begin{equation}
\label{eq:duaSDE}
dX(t) = f(t,X(t))dt + g(t,X(t))dB(t).
\end{equation}
Therefore, one can construct numerical methods for \eqref{eq:duatcSDE} by combining the numerical method to \eqref{eq:duaSDE} with the discretization of $E(t)$.
\par
However, such a one-to-one duality principle does not exist  for many other forms of time-changed SDEs. Thus, discretizing the time-changed SDE directly cannot be avoided for more general time-changed SDEs. For example, the time-changed SDEs of the form
\begin{equation}
\label{eq:noduatcSDE}
dY(t) = f(t,Y(t))dE(t) + g(t,Y(t))dB(E(t))
\end{equation}
do not have the dual classical SDEs. Therefore, the EM method was proposed in \cite{JK2019} by the direct discretization of the equation, where the assumption imposed on both $f$ and $g$ is the global Lipschitz condition in terms of the state variable. Such a requirement excludes many time-changed SDEs, like
\begin{equation}
\label{eq:nltcSDE}
dY(t) = (2Y(t)-9Y^5(t)) dE(t) + Y^3(t)d B(E(t)),
\end{equation}
where some super-linear terms appear in both $f$ and $g$. In the case of the classical SDE, i.e. $E(t)$ is replaced by $t$ in \eqref{eq:nltcSDE}, it was proved that the classical EM method fails to converge to such equations \cite{HJK2011}. There are many alternative methods have been proposed to handle the super-linearity, for example, the semi-implicit EM method \cite{HHKW2020,Hu1996,MS2013}, the tamed-type methods \cite{DFFM2021,HJK2012,Sabanis2013,SLL2018,WG2013}, the truncated-type methods \cite{CZG2020,LW2019,LMY2019,Mao2015,Zhang2020}, and the adaptive time-stepping methods \cite{FG2020,IJE2015,KL2018,NZB2015}. We just mention some of works here and refer the readers to the references therein.
\par
Although there is no existing work like \cite{HJK2011} that directly proves the divergence of the EM method for the time-changed SDE with super-linear growing coefficients, it is not hard to see by using the similar arguments in \cite{HJK2011} that the EM method is not convergent for this type of time-changed SDEs. To be slightly more precise, this is because that the time-changed Brownian motion is an unbounded stochastic process and would lead to the divergence of the EM method just like what the Brownian motion did for the EM method in the case of the classical SDEs. Therefore, some alternatives are needed to deal with time-changed SDEs like \eqref{eq:nltcSDE}.
\par
In this paper, we adopt the truncating idea to develop the truncated EM method for a class of highly non-linear time-changed SDEs. To our best knowledge, this paper is the first work to study numerical methods for highly non-linear time-changed SDEs, which do not have dual classical SDEs.
\par
The main contributions of this work are as follows.
\begin{itemize}
\item Strong convergence of the proposed method is proved when the state variables satisfy the polynomial growth condition and the time variables obey the H\"older continuity condition.
\item The convergence rate is obtained to be $\min\{\gamma_f,~\gamma_g,~\frac{1}{2} - \varepsilon \}$ for arbitrarily small $\varepsilon > 0$, where $\gamma_f$ and $\gamma_g$ are the H\"older index for the time variables.
\end{itemize}
\par
The rest of this paper is arranged in the following way. The mathematical preliminaries are put in Section 2. The main results and proofs are presented in Section 3. Numerical examples that illustrate the theoretical results are provided in Section 4. Section 5 concludes the paper and discusses some future research.

\section{Mathematical preliminaries}
This section is divided into two parts. The notations and assumptions are presented in Section 2.1. The truncated EM method is constructed in Section 2.2.
\subsection{Notations and assumptions}
The following setting-ups of probability spaces and the stochastic processes are quite standard. But to keep the paper self-contained, we decide to include them in this subsection.
\par
Throughout this paper, unless otherwise specified, we let $(\Omega_B,\F^B,\PP_B)$ be a complete probability space with a filtration $\{\F^B_t\}_{t \geq 0}$ satisfying the usual conditions (that is, it is right continuous and increasing while $\F^B_0$ contains all $\PP_B$-null sets). Let $B(t)= (B_1(t), B_2(t), ..., B_m(t))^T$ be an $m$-dimensional
$\F^B_t$-adapted standard Brownian motion. Let $\EE_B$ denote the probability expectation with respect to $\PP_B$. If $x\in \RR^d$, then $|x|$ is the Euclidean norm. Let $x^T$ denote the transposition of $x$. Moreover, for two real numbers $a$ and $b$, we use $a\vee b=\max(a,b)$ and $a\wedge b=\min(a,b)$.

\par
Let $D(t)$ be an $\F^D_t$-adapted subordinator (without killing), i.e. a nondecreasing L\'{e}vy process on $[0,\infty)$ starting
at $D(0)=0$  defined on a complete probability space $(\Omega_D , \F^D, \PP_D)$ with a filtration $\left\{\F^D_t\right\}_{t \ge 0}$ satisfying the usual conditions. Let $\EE_D$ denote the expectation under the probability measure $\PP_D$.

\par
The Laplace transform of $D(t)$ is of the form
$$
    \EE e^{-rD(t)} = e^{-t \phi(r)},\quad r>0,\,t\geq0,
$$
where the characteristic (Laplace) exponent $\phi:(0,\infty)\rightarrow(0,\infty)$ is a Bernstein function with $\phi(0+):=\lim_{r\downarrow0}\phi(r)=0$, i.e.\ a $C^\infty$-function such that $(-1)^{n-1}\phi^{(n)}\geq0$ for all $n \in \NN$.
Every such $\phi$ has a unique L\'{e}vy--Khintchine representation
$$
    \phi(r)
    =\vartheta r+\int_{(0,\infty)}\left(1-e^{-rx}\right) \nu(d x),\quad r>0,
$$
where $\vartheta\geq 0$ is the drift parameter and $\nu$ is a L\'{e}vy measure on $(0,\infty)$ satisfying $\int_{(0,\infty)}(1\wedge x) \nu(d x)<\infty$. We will focus on the case that $t\mapsto D(t)$ is a.s.\ strictly increasing, i.e.\ $\vartheta>0$ or $\nu(0,\infty)=\infty$; obviously, this is also equivalent to $\phi(\infty):=\lim_{r\rightarrow\infty}\phi(r)=\infty$.
\par
Let $E(t)$ be the (generalized, right-continuous) inverse of $D(t)$, i.e.
\begin{equation*}
E(t) := \inf\{ s\geq0\,;\,D(s) > t \}, \quad t \geq 0.
\end{equation*}
We call $E(t)$ an inverse subordinator associated with the Bernstein function $\phi$. Note that $t\mapsto E(t)$
is a.s.\ continuous and nondecreasing.
\par
We always assume that $B(t)$ and $D(t)$ are independent.  Define the product probability space by
\begin{equation*}
(\Omega , \F, \PP):= (\Omega_B \times \Omega_D, \F^B\otimes \F^D, \PP_B \otimes \PP_D).
\end{equation*}
Let $\EE$ denote the expectation under the probability measure $\PP$. It is clear that $\EE(\cdot) = \EE_D \EE_B (\cdot) = \EE_B \EE_D (\cdot)$.
The process $B(E(t))$ is called a time-changed Brownian motion and $B(E(t))$ is understood as a subdiffusion \cite{MS2004,UHK2018}.
\par
In this paper, we consider the d-dimension time-changed SDEs of the form
 \begin{equation}\label{SDE}
dY(t) = f(t,Y(t))dE(t) + g(t,Y(t))dB(E(t)),~~Y(0) = Y_0,~~t \geq 0,
\end{equation}
with $\EE |Y_0|^q < \infty$ for all $q > 0$. Here, $f: \RR_+ \times \RR^d \rightarrow \RR^d$ and $g: \RR_+ \times \RR^d \rightarrow \RR^{d\times m}$.
\par
The following assumptions are imposed on the coefficients of \eqref{SDE}.

\begin{assp}
\label{asspumption1}
Assume that there exist positive constants $\a$ and $L$ such that
\begin{equation*}
|f(t,x)-f(t,y)|\vee |g(t,x)-g(t,y)|\leq L(1+|x|^{\a}+|y|^{\a})|x-y|,
\end{equation*}
for all $t\in[0,T]$, any $x,y \in \RR^d$.
\end{assp}
\par
It can be observed from Assumption \ref{asspumption1} that all $t\in [0,T]$, and $x\in \RR^d$
\begin{equation}
\label{equ2-2}
|f(t,x)|\vee|g(t,x)|\leq M(1+|x|^{\a+1}),
\end{equation}
where $M$ depends on $L$ and $\sup_{0\leq t\leq T}\left(|f(t,0)|+|g(t,0)|\right)$.
\begin{assp}
\label{assumption1-2}
Assume that there exists a pair of constants $p>2$ and $K>0$ such that
\begin{equation*}
(x-y)^{\mathrm{T}}(f(t,x)-f(t,y)) +\frac{5p-1}{2}|g(t,x)-g(t,y)|^2\leq K|x-y|^2,
\end{equation*}
for all $t\in [0,T]$ and any $x,y \in \mathbb{R}^d$.
\end{assp}
\begin{assp}
\label{assumption2}
Assume that there exists a pair of constants $q>2$ and $K_1>0$ such that
\begin{equation*}
x^{\mathrm{T}}f(t,x) +\frac{5q-1}{2}|g(t,x)|^2\leq K_1(1+|x|^2),
\end{equation*}
for all $t\in [0,T]$ and any $x \in \mathbb{R}^d$.
\end{assp}
\begin{assp}
\label{assumption2-4}
Assume that there exists constants $\gamma_{f}\in (0,1]$, $\gamma_{g}\in (0,1]$, $H_1>0$ and $H_2>0$ such that
\begin{eqnarray*}
&&|f(s,x)-f(t,x)|\leq H_1(1+|x|^{\a+1})(s-t)^{\gamma_{f}},\\
&&|g(s,x)-g(t,x)|\leq H_2(1+|x|^{\a+1})(s-t)^{\gamma_{g}},
\end{eqnarray*}
for all $x,y \in \mathbb{R}^d$ and $s,t\in [0,T]$.
\end{assp}
\par
Under those assumptions above, the existence and uniqueness of the strong solution to \eqref{SDE} can be derived in the similar approach as that of Lemma 4.1 in \cite{Kobayashi2011}. It should be noted that the global Lipschitz condition is assumed in Lemma 4.1 of \cite{Kobayashi2011} though, its proof does not require this assumption explicitly. In fact, $Y(t)$ can be understood as a $\G_t$-semimartingale, where $\G_t = \F^B_{E(t)}$. Thus, the existence and uniqueness of the strong solution to \eqref{SDE} can be derived from  the existence and uniqueness of the strong solution to SDEs driven by semimartingale, see for example \cite{Mao1991} and \cite{Protter2004}.

\subsection{The truncated EM method for time-changed SDEs}
There are three steps to define the truncated EM method for time-changed SDEs.
\par \noindent
{\bf The first step} is to discretise the inverse subordinator $E(t)$ in a time interval $[0,T]$ for any given $T>0$. Following the idea in \cite{GM2010}, we simulate the path of $D(t)$ by $D_\Delta(t_i) = D_\Delta(t_{i-1} )+ h_i$ with $D_\Delta(0) = 0$ firstly, where $h_i$ is independently identically sequence with $h_i = D(\Delta)$ in distribution and $\Delta>0$ is a constant. The procedure is stopped when
\begin{equation*}
T \in [ D_\Delta(t_{N}), D_\Delta(t_{N+1})),
\end{equation*}
for some $N$. Then the approximate $E_\Delta(t)$ to $E(t)$ is generated by
\begin{equation}\label{findEht}
E_\Delta(t) = \big(\min\{n; D_\Delta(t_n) > t\} - 1\big)\Delta,
\end{equation}
for $t \in [0,T]$. It is easy to see
\begin{equation*}
E_\Delta(t) = i\Delta,\quad\text{when $t \in   \left[ D_\Delta(t_{i}), D_\Delta(t_{i+1})\right)$}.
\end{equation*}
For $i = 1,2, ..., N$, denote
\begin{equation*}
\rho_i = D_\Delta(t_i).
\end{equation*}
It is clear that
\begin{equation}\label{eq:EDelta}
E_\Delta(\rho_i) = E_\Delta(D_h(t_i)) = i\Delta.
\end{equation}
\par \noindent
{\bf The second step} is to construct the truncating function. We choose a strictly increasing continuous function $\m:\RR_+\rightarrow \RR_+$ such that $\m(u)\rightarrow \infty$ as $u\rightarrow \infty$ and
\begin{equation*}
\sup_{0\leq t\leq T}\sup_{|x|\leq u}(|f(t,x)|\vee|g(t,x)|)\leq \m(u),\quad \forall u\geq 1.
\end{equation*}
Denote by $\m^{-1}$ the inverse function of $\m$. It is clear that $\m^{-1}$ is a strictly increasing continuous function from $[\m(0),\infty)$ to $\RR_+$. We also choose a constant $\hat{\k}\geq 1\wedge \m(1)$ and a strictly decreasing function $\k:(0,1]\rightarrow [\m(1),\infty)$ such that
\begin{equation}
\label{equ0}
\lim_{\D \rightarrow 0}\k(\D)=\infty,\quad \D^{1/4}\k(\D)\leq \hat{\k},\quad \forall \D \in (0,1].
\end{equation}
\par
For a given step size $\D \in (0,1]$, define the truncated mapping $\pi_{\D}:\RR^d \rightarrow \{ x\in \RR^d:|x|\leq \m^{-1}(\k(\D))\} $ by
\begin{equation*}
\pi_{\D}(x)=\left(|x|\wedge\m^{-1}(\k(\D))\right)\frac{x}{|x|},
\end{equation*}
where we set $x/|x|=0$ when $x=0$.
Define the truncated functions by
\begin{equation*}
f_{\D}(t,x)=f(t,\pi_{\D}(x)),\quad g_{\D}(t,x)=g(t,\pi_{\D}(x))
\end{equation*}
for $x\in \RR^d$. It is easy to see that for any $t\in[0,T]$ and all $x \in \RR^d$,
\begin{equation}
\label{equ01}
|f_{\D}(t,x)|\vee |g_{\D}(t,x)|\leq \m(\m^{-1}(\k(\D)))=\k(\D).
\end{equation}
\par \noindent
{\bf The last step} is to discretise the equation. We define the discrete version of the  numerical method at the time grid $\rho_n$ to get
\begin{eqnarray*}
X_{\r_{n+1}}&=&X_{\r_{n}}+f_{\D}(\r_n,X_{\r_{n}})\bigg(E_\Delta(\rho_{n+1}) - E_\Delta(\rho_n)\bigg)\nonumber \\
&&+g_{\D}(\r_n,X_{\r_{n}})\bigg(B(E_\Delta(\rho_{n+1})) - B(E_\Delta(\rho_n)) \bigg).
\end{eqnarray*}
It should be noted that $\{\r_n\}_{n=1,2,...,N}$ is a random sequence but independent from the Brownian motion. In addition, It is not hard to see from \eqref{eq:EDelta} that
\begin{equation*}
E_\Delta(\rho_{n+1}) - E_\Delta(\rho_n) = \Delta
\end{equation*}
and
\begin{equation*}
B(E_\Delta(\rho_{n+1})) - B(E_\Delta(\rho_n)) = B((n+1)\Delta) -B(n\Delta).
\end{equation*}
Therefore, the discrete version of the truncated EM method is defined a discrete-time process $(X_{\r_n})_{n\in \{0,1,2,\cdots,N\}}$ by setting $X_0=Y(0)$ and
\begin{equation}\label{eq:DiscreteEM}
X_{\r_{n+1}}=X_{\r_{n}}+f_{\D}(\r_n,X_{\r_{n}})\D+g_{\D}(\r_n,X_{\r_{n}})[B((n+1)\D)-B(n\D)].
\end{equation}
\par
For the convenience of the analysis, two versions of the continuous-time truncated EM methods are formed. The first one is defined by
\begin{equation}
\label{equ03}
\bar{X}(t)=\sum_{n=0}^{N\D}X_{\r_n}\ind_{[\r_n,\r_{n+1})}(t),
\end{equation}
which is a simple step process. The second one is defined by
\begin{equation}
\label{equ04}
X(t)=X(0)+\int_{0}^{t}f_{\D}{(\bar{\r}(s),\bar{X}(s))dE(s)}+\int_{0}^{t}g_{\D}{(\bar{\r}(s),\bar{X}(s))dB(E(s))},
\end{equation}
which is continuous in $t\in[0, T]$, where $X(0)=Y(0)$ and $\bar{\r}(s)=\r_n\ind_{[\r_n,\r_{n+1})}(s)$.

\section{Main results}
This section consists of three parts. Our main theorem and an important corollary are presented in Section 3.1, while the proof of the main theorem is postponed to Section 3.3. Several useful lemmas are proved in Section 3.2 as the preparations for the proof of the main theorem.

\subsection{Main theorem and its corollary}
\begin{theorem}
\label{theorem3-2}
Let Assumptions \ref{asspumption1}, \ref{assumption1-2}, \ref{assumption2} and \ref{assumption2-4} hold with $q>(\a+1)p$. Then, for any $\bar{p}\in [2,p)$ and $\D \in (0,1]$
\begin{equation}
\label{th321}
\EE \bigg( \sup_{0\leq t\leq T}|Y(t)-X(t)|^{\bar{p}} \bigg)\leq C\bigg(\D^{\gamma_f\bar{p}}+\D^{\gamma_g\bar{p}}+\D^{\bar{p}/2}(\k(\D))^{\bar{p}}+(\m^{-1}(\k(\D)))^{(\a+1)\bar{p}-q}\bigg),
\end{equation}
and
\begin{equation}
\label{th322}
\EE\bigg(\sup_{0\leq t\leq T}|Y(t)-\bar{X}(t)|^{\bar{p}}\bigg)\leq C\bigg(\D^{\gamma_f\bar{p}}+\D^{\gamma_g\bar{p}}+\D^{\bar{p}/2}(\k(\D))^{\bar{p}}+(\m^{-1}(\k(\D)))^{(\a+1)\bar{p}-q}\bigg).
\end{equation}
\end{theorem}
By strengthening Assumption \ref{assumption2} and choosing specific $\mu(\cdot)$ and $\kappa (\cdot)$, we achieve the following corollary, which displays the strong convergence rate more clearly.

\begin{coro}
\label{theorem3-1}
Let Assumptions \ref{asspumption1}, \ref{assumption1-2} and \ref{assumption2-4} hold, and let Assumption \ref{assumption2} hold for any $q>2$. In particular, recalling \eqref{equ2-2}, we may define
\begin{equation}
\label{mu}
\m(u)=2Mu^{\a+2},u\geq 1,
\end{equation}
and let
\begin{equation}
\label{k}
 \k(\D)=\D^{-\varepsilon} \quad \text{for some}\quad \varepsilon\in(0,\frac{1}{4}]\quad  \text{and}\quad \hat{\k}\geq 1.
\end{equation}
Then for any $\bar{p}\geq 2$, $\D \in (0,1]$ and $\varepsilon \in (0,\frac{1}{4}]$,
\begin{equation}
\label{th311}
\EE\left( \sup_{0\leq t\leq T}|Y(t)-X(t)|^{\bar{p}}\right)\leq C\D^{\min\{\gamma_f\bar{p},\gamma_g\bar{p},(\frac{1}{2}-\varepsilon)\bar{p}\}}
\end{equation}
and
\begin{equation}
\label{th312}
\EE\left( \sup_{0\leq t\leq T}|Y(t)-\bar{X}(t)|^{\bar{p}}\right)\leq C\D^{\min\{\gamma_f\bar{p},\gamma_g\bar{p},(\frac{1}{2}-\varepsilon)\bar{p}\}}.
\end{equation}
\end{coro}
\begin{proof}
By \eqref{mu}, we can get
\begin{equation*}
\m^{-1}(u)=\left(\frac{u}{2M}\right)^{\frac{1}{\a+2}}.
\end{equation*}
It is derived from Theorem \ref{theorem3-2} that
\begin{equation*}
\EE\left(\sup_{0\leq t\leq T}|Y(t)-X(t)|^{\bar{p}}\right)\leq C\bigg(\D^{\gamma_f\bar{p}}+\D^{\gamma_g\bar{p}}+\D^{\frac{(1-2\varepsilon)\bar{p}}{2}}+\D^{\frac{\varepsilon(q-(\a+1)\bar{p})}{\a+2}}\bigg),
\end{equation*}
and
\begin{equation*}
\EE\left(\sup_{0\leq t\leq T}|Y(t)-\bar{X}(t)|^{\bar{p}}\right)\leq C\bigg(\D^{\gamma_f\bar{p}}+\D^{\gamma_g\bar{p}}+\D^{\frac{(1-2\varepsilon)\bar{p}}{2}}+\D^{\frac{\varepsilon(q-(\a+1)\bar{p})}{\a+2}}\bigg).
\end{equation*}
Then, choosing $q$ sufficiently large such that
\begin{equation*}
\frac{\varepsilon(q-(\a+1)\bar{p})}{\a+2}>\min\{\gamma_f\bar{p},\gamma_g\bar{p},(\frac{1}{2}-\varepsilon)\bar{p}\},
\end{equation*}
we can draw the assertions \eqref{th311} and \eqref{th312} immediately. \eproof
\end{proof}

\begin{rmk}
Corollary  \ref{theorem3-1} states that the strong convergence rate of the truncated EM method can be arbitrarily close to 0.5 for some smooth coefficients in terms of the time variable. However, with the coefficients becoming more singular the convergence rate gets lower.
\end{rmk}

\subsection{Useful lemmas}
Five lemmas are needed before we proceed to the proof of the main theorem.
\par
The first two lemmas in this subsection state that the truncated functions $f_\Delta$ and $g_\Delta$ inherit Assumptions  \ref{asspumption1} and \ref{assumption2} to some extended. Since the proofs of these two lemmas are similar to those of Lemmas 3.2 and 3.5 in \cite{GLM2018}, we omit them here.
\begin{lemma}
\label{lemma23}
Let Assumption \ref{assumption2} hold. Then, for all $\D\in (0,1]$, we have
\begin{equation*}
 x^{\mathrm{T}}f_{\D}(t,x) +\frac{5q-1}{2}|g_{\D}(t,x)|^2 \leq \hat{K}_1(1+|x|^2),\quad \forall x\in \RR^d,
\end{equation*}
where $\hat{K}_1=2K_1\left(1\vee \frac{1}{\m^{-1}(\k(1))}\right)$.
\end{lemma}

\begin{lemma}
\label{lemma2-6}
Let Assumption \ref{asspumption1} hold. Then, for all $\D\in (0,1]$, we have
\begin{equation*}
|f_{\D}(t,x)-f_{\D}(t,y)|\vee |g_{\D}(t,x)-g_{\D}(t,y)|\leq L(1+|x|^{\a}+|y|^{\a})|x-y|,
\end{equation*}
for all $t\in(0,T]$, $x,y \in \RR^{d}$.
\end{lemma}

The lemma below states the moment boundedness of the underlying solution.
\par
\begin{lemma}
\label{lemY}
Suppose Assumption \ref{asspumption1} and \ref{assumption2} hold. Then, for any $p\in [2,q)$
\begin{equation*}
\EE \left( \sup _{0\leq t\leq T}|Y(t)|^p \right)<\infty.
\end{equation*}
\end{lemma}
\par
\begin{proof}
Denote $\tau_\ell :=\inf\{t\geq0;|Y(t)|>\ell \}$ for some positive integer $\ell$. It can be seen that
\begin{equation*}
\int_{0}^{t}{\EE_{B}\left(\sup_{0\leq s \leq r\wedge \tau_\ell}|Y(s)|^p\right) dE(r)< \ell^p E(t)< \infty }.
\end{equation*}
\par
Fix $\ell$ and $t\in [0,T]$. Since E has continuous paths, $B\circ E$ is a continuous martingale with quadratic variation $[B\circ E,B\circ E]=[B,B]\circ E=E$. By the It\^o formula, we can derive from \eqref{SDE} that
\begin{equation}
\label{equ}
|Y(s)|^p=|Y(0)|^p + A_0+ M_s,
\end{equation}
where
\begin{equation*}
A_0:= \int_{0}^{s}{\left\{p|Y(r)|^{p-1}f(r,Y(r))+\frac{1}{2}p(p-1)|Y(r)|^{p-2}|g(r,Y(r))|^2\right\}dE(r)}
\end{equation*}
and
\begin{equation*}
M_s:=\int_{0}^{s}{p|Y(r)|^{p-1}g(r,Y(r))dB(E(r))}.
\end{equation*}
It can be noted that the stochastic integral $(M_t)_{t\geq 0}$ is a local martingale with quadratic variation
\begin{equation*}
[M,M]_t=\int_{0}^{t}{p^2|Y(r)|^{2p-2}|g(r,Y(r))|^2dE(r)}.
\end{equation*}
For $0\leq r \leq t\wedge \tau_\ell$, we have
\begin{equation*}
\begin{split}
p^2|Y(r)|^{2p-2}|g(r,Y(r))|^2&=p^2|Y(r)|^{p}|Y(r)|^{p-2}|g(r,Y(r))|^2\\
&\leq p^2\left(1+\sup_{0\leq s\leq t\wedge \tau_\ell}|Y(s)|^p\right)|Y(r)|^{p-2}|g(r,Y(r))|^2.
\end{split}
\end{equation*}
By using the inequality $(ab)^{1/2}\leq a/\l+\l b$ valid for any $a,b>0$ and $\l >0$, we can see that for $0 < s < t\wedge \tau_\ell$
\begin{equation}
\label{equ2}
\begin{split}
&([M,M]_{s})^{1/2}\\
\leq& p\left((1+\sup_{0\leq s\leq t\wedge \tau_\ell}|Y(s)|^p)\int_{0}^{s}{|Y(r)|^{p-2}|g(r,Y(r))|^2 dE(r)}\right)^{\frac{1}{2}}\\
\leq& p\left(\frac{(1+\sup_{0\leq s\leq t\wedge \tau_\ell}|Y(s)|^p)}{2pb_1}+2pb_1\int_{0}^{s}{|Y(r)|^{p-2}|g(r,Y(r))|^2 dE(r)}\right),
\end{split}
\end{equation}
where $b_1$ is the constant appearing in the Burkholder--Davis--Gundy (BDG) equality (Chapter 1, Theorem 7.3 in \cite{Mao2008}). Substituting \eqref{equ2} into \eqref{equ} and taking $\EE_B$ on both sides, we have
\begin{equation*}
\begin{split}
&1+\EE_B\left(\sup_{0\leq s\leq t\wedge \tau_\ell}|Y(s)|^p\right)\\
\leq& 1+|Y(0)|^p+\EE_B\left(\sup_{0\leq s\leq t\wedge \tau_\ell}\int_{0}^{s}{\{p|Y(r)|^{p-2}Y^T(r)f(r,Y(r))+\frac{1}{2}p(p-1)|Y(r)|^{p-2}|g(r,Y(r))|^2\}dE(r)}\right)\\
\quad&+b_1\EE_B\left(\int_{0}^{s}{2p^{2}b_1|Y(r)|^{p-2}|g(r,Y(r))|^2 dE(r)}+\frac{1}{2b_1}(1+\sup_{0\leq s\leq t\wedge \tau_\ell}|Y(s)|^p)\right)\\
=&1+|Y(0)|^p+\frac{1}{2}\left(1+\EE_B(\sup_{0\leq s\leq t\wedge \tau_\ell}|Y(s)|^p)\right)\\
\quad&+b\EE_B\left(\sup_{0\leq s\leq t\wedge \tau_\ell}\int_{0}^{s}{p|Y(r)|^{p-2}\left(Y^{\mathrm{T}}(r)f(r,Y(r)) +\frac{5p-1}{2}|g(r,Y(r))|^2 \right)dE(r)}\right),
\end{split}
\end{equation*}
where $b=b_1\vee 1$. Noting that for any nonnegative process $a(r)$, the inequality
\begin{equation}
\label{stti}
\int_{0}^{t\wedge \tau_\ell}{a(r)dE(r)}\leq \int_{0}^{t}{a(r\wedge \tau_\ell)dE(r)}
\end{equation}
holds. Using Assumption \ref{assumption2}, we obtain
\begin{equation*}
\begin{split}
1+\EE_B\left(\sup_{0\leq s\leq t\wedge \tau_\ell}|Y(s)|^p\right)&\leq 1+|Y(0)|^p+\frac{1}{2}\left(1+\EE_B(\sup_{0\leq s\leq t\wedge \tau_\ell}|Y(s)|^p)\right)\\
&\quad+b\EE_B\left(\int_{0}^{t}{pK_1|Y(r\wedge \tau_\ell)|^{p-2}\left(1+|Y(r)|^2\right)dE(r)}\right)\\
&\leq 1+|Y(0)|^p+\frac{1}{2}\left(1+\EE_B(\sup_{0\leq s\leq t\wedge \tau_\ell}|Y(s)|^p)\right)\\
&\quad+2pbK_1\int_{0}^{t}{\left(1+\EE_B(\sup_{0\leq s\leq r\wedge \tau_\ell}|Y(s)|^p)\right)dE(r)},
\end{split}
\end{equation*}
which in turn gives
\begin{equation*}
1+\EE_B \left(\sup_{0\leq s\leq t\wedge \tau_\ell}|Y(s)|^p\right)\leq 2(1+|Y(0)|^p)+4pbK_1\int_{0}^{t}\left(1+\EE_B\left(\sup_{0\leq s\leq r\wedge \tau_\ell}|Y(s)|^p\right)\right)dE(r).
\end{equation*}
\par
Thus, applying the Gronwall-type inequality (Chapter IX.6a, Lemma 6.3 in \cite{JS2003}) yields
\begin{equation*}
1+\EE_B \left(\sup_{0\leq s\leq t\wedge \tau_\ell}|Y(s)|^p\right)\leq 2(1+|Y(0)|^p)e^{4pbK_1E(t)}
\end{equation*}
for all $t\in [0,T]$.
\par
Since $\tau_{\ell} \rightarrow \infty$ as $\ell \rightarrow \infty$, setting $t=T$ and letting $\ell\rightarrow \infty$ give
\begin{equation*}
1+\EE_B\left(\sup_{0\leq t\leq T}|Y(s)|^p\right)\leq 2(1+|Y(0)|^p)e^{4pbK_1E(T)}.
\end{equation*}
Taking $\EE_D$ on both sides and using the fact that $\EE_D [e^{\l E(T)}]<\infty$ for any $\l >0$ ( \cite{DL2020online}, \cite{JK2016}, \cite{MOW2011} ) yield the desired result.
\eproof
\end{proof}
The next lemma measures the difference between the two versions of the continuous-time truncated EM method.
\begin{lemma}
\label{lemma2}
For any $\D\in (0,1]$ and any $\hat{p}>2$, we have
\begin{equation}
\label{equ06}
\EE_B|X(t)-\bar{X}(t)|^{\hat{p}}\leq c_{\hat{p}}\D^{\hat{p}/2}\left(\k(\D)\right)^{\hat{p}},\quad \forall t\geq 0,
\end{equation}
where $c_{\hat{p}}=b_{\hat{p}}2^{\hat{p}-1}$. Consequently,
\begin{equation}
\label{equ07}
\lim_{\D\rightarrow 0}\EE_B|X(t)-\bar{X}(t)|^{\hat{p}}=0,\quad \forall t\geq 0.
\end{equation}
\end{lemma}
\begin{proof}
Fix any $\D\in (0,1]$, $\hat{p}>2$ and $t\geq 0$. There is a unique integer $n\geq 0$ such that $\r_n\leq t\leq \r_{n+1}$. By properties of the It\^o integral, we then derive from \eqref{equ04} that
\begin{equation}
\label{lem250}
\begin{split}
&\EE_B|X(t)-\bar{X}(t)|^{\hat{p}}\\
=& \EE_B|X(t)-X(\r_n)|^{\hat{p}}\\
\leq& 2^{\hat{p}-1}\left(\EE_B\left|\int_{\r_n}^{t}{f_{\D}(\bar{\r}(s),\bar{X}(s))dE(s)} \right|^{\hat{p}}+\EE_B\left|\int_{\r_n}^{t}{g_{\D}(\bar{\r}(s),\bar{X}(s))dB(E(s))} \right|^{\hat{p}}\right).\\
\end{split}
\end{equation}
The first item in the above brackets, using the H\"older inequality, we can see
\begin{equation}
\label{lem25}
\EE_B\left|\int_{\r_n}^{t}{f_{\D}(\bar{\r}(s),\bar{X}(s))dE(s)} \right|^{\hat{p}}\leq \D^{\hat{p}-1}\EE_B\int_{\r_n}^{t}{\left|f_{\D}(\bar{\r}(s),\bar{X}(s))\right|^{\hat{p}} dE(s)},
\end{equation}
the second item in the above brackets, using the Burkholder-Davis-Gundy inequality, we can see
\begin{equation}
\label{lem251}
\begin{split}
\EE_B\left|\int_{\r_n}^{t}{g_{\D}(\bar{\r}(s),\bar{X}(s))dB(E(s))}\right|^{\hat{p}} &=\EE_B\left(\left|\int_{\r_n}^{t}g_{\D}(\bar{\r}(s),\bar{X}(s))dB(E(s))\right|^2\right)^{\hat{p}/2}\\
&\leq b_{\hat{p}}\EE_B\left(\int_{\r_n}^{t}{\left|g_{\D}(\bar{\r}(s),\bar{X}(s))\right|^2 dE(s)}\right)^{\hat{p}/2}\\
&\leq b_{\hat{p}}\D^{\hat{p}/2-1}\EE_B\int_{\r_n}^{t}{\left|g_{\D}(\bar{\r}(s),\bar{X}(s))\right|^{\hat{p}} dE(s)}.
\end{split}
\end{equation}
Substituting the estimates \eqref{lem25} and \eqref{lem251} into \eqref{lem250}, we have
\begin{equation*}
\begin{split}
\EE_B|X(t)-\bar{X}(t)|^{\hat{p}} &\leq c_{\hat{p}}\left(\D^{\hat{p}-1}\D(\k(\D))^{\hat{p}}+\D^{\hat{p}/2-1}\D(\k(\D))^{\hat{p}}\right)\\
&\leq c_{\hat{p}}\D^{\hat{p}/2}(\k(\D))^{\hat{p}},
\end{split}
\end{equation*}
where \eqref{equ01} is used and $c_{\hat{p}}=b_{\hat{p}}2^{\hat{p}-1}$. This completes the proof of \eqref{equ06}. Noting from \eqref{equ0}, we have $\D^{\hat{p}/2}(\k(\D))^{\hat{p}}\leq \D^{\hat{p}/4}$. Then, \eqref{equ07} can be derived from \eqref{equ06}.
\eproof
\end{proof}
\par
The following lemma shows the moment boundedness of the truncated EM method.
\begin{lemma}
\label{lemma3}
Let Assumptions \ref{asspumption1} and \ref{assumption2} hold. Then
\begin{equation}
\label{equ08}
\sup_{0< \D \leq 1}\EE [\sup_{0\leq t\leq T}|X(t)|^p]\leq C, \quad \forall T>0,
\end{equation}
where $C$ is a constant dependent on $X(0)$, $p$, $T$ and $\hat{K}_1$, but independent from $\D$.
\end{lemma}
\begin{proof}
Define the stopping time $\zeta_{\ell}:=\inf\{ t\geq 0; |X(t)|>\ell \}$ for some positive integer $\ell$. It can be seen that
\begin{equation*}
\int_{0}^{t}{\EE_B\left(\sup_{0\leq s \leq t\wedge \zeta_{\ell}}|X(s)|^p\right)dE(r)}\leq \ell^pE(t)<\infty.
\end{equation*}
Fix any $\D\in(0,1]$ and $T\geq 0$. By the It\^o formula, we derive from \eqref{equ04} that, for $0\leq u\leq t\wedge \zeta_{\ell}$,
\begin{equation*}
|X(u)|^p=|X(0)|^p+A_1+M_u,
\end{equation*}
where
\begin{equation*}
A_1:=\int_{0}^{u}{\left(p|X(s)|^{p-2}X^{\mathrm{T}}(s)f_{\D}(\bar{\r}(s),\bar{X}(s))+\frac{1}{2}p(p-1)|X(s)|^{p-2}|g_{\D}(\bar{\r}(s),\bar{X}(s))|^2\right)dE(s)},
\end{equation*}
\begin{equation*}
M_u:=\int_{0}^{u}{p|X(s)|^{p-1}|g_{\D}(\bar{\r}(s),\bar{X}(s))|dB(E(s))}.
\end{equation*}
It can be noted that the stochastic integral $(M_t)_{t\geq 0}$ is a local martingale with quadratic variation
\begin{equation*}
[M,M]_t=\int_{0}^{t}{p^2|X(s)|^{2p-2}|g_{\D}(\bar{\r}(s),\bar{X}(s))|^2 dE(s)}.
\end{equation*}
For $0\leq s\leq t\wedge \zeta_{\ell}$,
\begin{equation*}
\begin{split}
p^2|X(s)|^{2p-2}|g_{\D}(\bar{\r}(s),\bar{X}(s))|^2&\leq p^2|X(s)|^p|X(s)|^{p-2}|g_{\D}(\bar{\r}(s),\bar{X}(s))|^2\\
&\leq p^2\left(\sup_{0\leq u\leq t\wedge \zeta_{\ell}}|X(u)|^p\right)|X(s)|^{p-2}|g_{\D}(\bar{\r}(s),\bar{X}(s))|^2.
\end{split}
\end{equation*}
Similar to the proof of Lemma \ref{lemY}. By using the inequality $(ab)^{1/2}\leq a/\l +\l b$ valid for any $a,b\geq 0$ and $\l \leq 0$, we can see that for $0\leq u \leq t\wedge \zeta_{\ell}$,
\begin{equation*}
\begin{split}
&([M,M]_t)^{1/2}\\
\leq& p\left(\sup_{0\leq u\leq t\wedge \zeta_{\ell}}|X(u)|^p\int_{0}^{u}{|X(s)|^{p-2}|g_{\D}(\bar{\r}(s),\bar{X}(s))|^2 dE(s)}\right)^{1/2}\\
\leq& p\left(\frac{\sup_{0\leq u\leq t\wedge \zeta_{\ell}}|X(u)|^p}{2pb_1}+2pb_1\int_{0}^{u}{|X(s)|^{p-2}|g_{\D}(\bar{\r}(s),\bar{X}(s))|^2 dE(s)}\right).
\end{split}
\end{equation*}
Therefore, for any $0\leq u\leq t\wedge \zeta_{\ell}$, by Lemma \ref{lemma23} and the Young inequality
\begin{equation*}
a^{p-2}b\leq \frac{p-2}{p}a^p + \frac{2}{p}b^{p/2},\quad \forall a,b \geq 0,
\end{equation*}
we then have
\begin{equation*}
\begin{split}
&\EE_B\left(\sup_{0\leq u\leq t\wedge \zeta_{\ell}}|X(u)|^p\right)\\
\leq& |X(0)|^p+\EE_B\bigg(\sup_{0\leq u\leq t\wedge \zeta_{\ell}}\int_{0}^{u}\bigg(p|X(s)|^{p-2}X^{\mathrm{T}}(s)f_{\D}(\bar{\r}(s),\bar{X}(s))\\
\quad& +\frac{1}{2}p(p-1)|X(s)|^{p-2}|g_{\D}(\bar{\r}(s),\bar{X}(s))|^2\bigg)\bigg)  dE(s)\\
\quad& +b_1\EE_B\left(\frac{1}{2b_1}\sup_{0\leq u\leq t\wedge \zeta_{\ell}}|X(u)|^p+\int_{0}^{u}{2b_1p^2|X(s)|^{p-2}|g_{\D}(\bar{\r}(s),\bar{X}(s))|^2 dE(s)}\right)\\
\leq& |X(0)|^p+\frac{1}{2}\EE_B\left(\sup_{0\leq u\leq t\wedge \zeta_{\ell}}|X(u)|^p\right)\\
\quad& +b\EE_B\left(\sup_{0\leq u\leq t\wedge \zeta_{\ell}}\int_{0}^{u}{p|X(s)|^{p-2}\left(\bar{X}^{\mathrm{T}}(s)f_{\D}(\bar{\r}(s),\bar{X}(s))+\frac{5p-1}{2}|g_{\D}(\bar{\r}(s),\bar{X}(s))|^2\right)dE(s)}\right)\\
\quad& +\EE_B\left(\int_{0}^{u}{p|X(s)|^{p-2}(X(s)-\bar{X}(s))^{\mathrm{T}}f_{\D}(\bar{\r}(s),\bar{X}(s))dE(s)}\right)\\
\leq& |X(0)|^p+\frac{1}{2}\EE_B\left(\sup_{0\leq u\leq t\wedge \zeta_{\ell}}|X(u)|^p\right)+bp\hat{K}_1\EE_B\left(\sup_{0\leq u\leq t\wedge \zeta_{\ell}}\int_{0}^{u}{|X(s)|^{p-2}(1+|\bar{X}(s)|^2)dE(s)}\right)\\
\quad& +\EE_B\left((p-2)\int_{0}^{u}{|X(s)|^p dE(s)}+2\int_{0}^{u}{|X(s)-\bar{X}(s)|^{p/2}|f_{\D}(\bar{\r}(s),\bar{X}(s))|^{p/2}dE(s)}\right),
\end{split}
\end{equation*}
where $b=b_1^2\vee 1$. Thus, for any $0\leq u\leq t\wedge \zeta_{\ell}$ and apply \eqref{stti} inequality, we have
\begin{equation*}
\begin{split}
\EE_B\left(\sup_{0\leq u\leq t\wedge \zeta_{\ell}}|X(t)|^p\right)&\leq 2|X(0)|^p+2bp\hat{K}_1\EE_B\int_{0}^{t}{|X(t\wedge \zeta_{\ell})|^{p-2}(1+|\bar{X}(s)|^2)dE(s)}\\
&\quad +2(p-2)\int_{0}^{t}{\EE_B |X(t\wedge \zeta_{\ell})|^p dE(s)}\\
&\quad +4\EE_B\int_{0}^{t}{|X(s)-\bar{X}(s)|^{p/2}|f_{\D}(\bar{\r}(s),\bar{X}(s))|^{p/2}dE(s)}.
\end{split}
\end{equation*}
By Lemma \ref{lemma2}, inequalities \eqref{equ01} and \eqref{equ0}, we have
\begin{equation*}
\label{}
\begin{split}
&\EE_B\int_{0}^{t}{|X(s)-\bar{X}(s)|^{p/2}|f_{\D}(\bar{\r}(s),\bar{X}(s))|^{p/2}dE(s)}\\
\leq& \left(\k(\D)\right)^{p/2}\int_{0}^{t}{\EE_B\left(|X(s)-\bar{X}(s)|^{p/2}\right)dE(s)}\\
\leq& \left(\k(\D)\right)^{p/2}\int_{0}^{t}{\left(\EE_B|X(s)-\bar{X}(s)|^p\right)^{1/2}dE(s)}\\
\leq& 2^{p-1}E(t)\left(\k(\D)\right)^p\D^{p/4}\\
\leq& 2^{p-1}E(t).
\end{split}
\end{equation*}
We therefore have
\begin{equation*}
\begin{split}
\EE_B\left(\sup_{0\leq u\leq t\wedge \zeta_{\ell}}|X(u)|^p\right)&\leq C_1+C_2\int_{0}^{t}{\left(\EE_B|X(t\wedge \zeta_{\ell})|^p+\EE_B|\bar{X}(t\wedge \zeta_{\ell})|^p\right)dE(s)}\\
&\leq C_1+2C_2\int_{0}^{t}{\EE_B\left(\sup_{0\leq u\leq t\wedge \zeta_{\ell}}|X(u)|^p\right) dE(s)},
\end{split}
\end{equation*}
where $C_1=2|X(0)|^p+2^{p+1}E(t)$ and $C_2=2bp\hat{K}_1\vee 2(p-2)$. Applying the well-known Gronwall-type inequality in Chapter IX.6a, Lemma 6.3 of \cite{JS2003} yields, for any $t\in[0,T]$,
\begin{equation*}
\EE_B\left(\sup_{0\leq u\leq t\wedge \zeta_{\ell}}|X(u)|^p\right)\leq C_1e^{2C_2E(t)}.
\end{equation*}
Since $\zeta_{\ell}\rightarrow \infty$ as $\ell\rightarrow \infty$. Setting $t=T$ and letting $\ell\rightarrow \infty$ give
\begin{equation*}
\EE_B\left(\sup_{0\leq t\leq T}|X(t)|^p\right)\leq C_1e^{2C_2E(T)}.
\end{equation*}
Taking $\EE_D$ on both sides, and using the fact that $\EE_D\left(E(T)e^{E(T)}\right)<\EE_D\left(e^{2 E(T)}\right)< \infty$ yield,
\begin{equation*}
\EE\left(\sup_{0\leq t\leq T}|X(t)|^p\right) \leq C.
\end{equation*}
As this holds for any $\D\in(0,1]$ and $C$ is independent of $\D$, we see the required assertion \eqref{equ08}. \eproof
\end{proof}
Now, we are ready to prove the main theorem.
\subsection{Proof of Theorem \ref{theorem3-2}}
\noindent
{\bf Proof.} Fix $\bar{p}\in [1,p)$ and $\D\in(0,1]$ arbitrarily. Let $e_{\D}(t)=Y(t)-X(t)$ for $t\geq 0$. For each integer $\ell> |Y(0)|$, define the stopping time
\begin{equation}
\theta_{\ell}=\inf\{t\geq0:|Y(t)|\vee|X(t)|\geq \ell\},
\end{equation}
where we set $\inf\emptyset=\infty$ (as usual $\emptyset$ denotes the empty set). By the It\^o formula, we have that for any $0\leq t\leq T$,
\begin{equation}
\label{error}
\begin{split}
|e_{\D}(t\wedge \theta_{\ell})|^{\bar{p}}&=\int^{t\wedge\theta_{\ell}}_{0}\bigg(\bar{p}|e_{\D}(s)|^{\bar{p}-1}\left(f(s,Y(s))-f_{\D}(\bar{\rho}(s),\bar{X}(s))\right)\\
&\quad +\frac{\bar{p}(\bar{p}-1)}{2}|e_{\D}(s)|^{\bar{p}-2}|g(s,Y(s))-g_{\D}(\bar{\rho}(s),\bar{X}(s))|^2 \bigg)dE(s)+M_{t\wedge \theta_{\ell}},
\end{split}
\end{equation}
where
\begin{equation*}
M_{t\wedge \theta_{\ell}}:=\int^{t\wedge\theta_{\ell}}_{0}\bar{p}|e_{\D}(s)|^{\bar{p}-1}|g(s,Y(s))-g_{\D}(\bar{\rho}(s),\bar{X}(s))|dB(E(s)).
\end{equation*}
Note that the stochastic integral $(M_t)_{t\geq 0}$ is a local martingale with quadratic variation
\begin{equation*}
[M,M]_{t\wedge \theta_{\ell}}=\int^{t\wedge\theta_{\ell}}_{0}\bar{p}^2|e_{\D}(s)|^{2\bar{p}-2}|g(s,Y(s))-g_{\D}(\bar{\rho}(s),\bar{X}(s))|^2 dE(s).
\end{equation*}
For $0\leq s\leq t\wedge\theta_{\ell}$, we have
\begin{equation*}
\begin{split}
&\bar{p}^2|e_{\D}(s)|^{2\bar{p}-2}|g(s,Y(s))-g_{\D}(\bar{\rho}(s),\bar{X}(s))|^2\\
=&\bar{p}^2|e_{\D}(s)|^{\bar{p}}|e_{\D}(s)|^{\bar{p}-2}|g(s,Y(s))-g_{\D}(\bar{\rho}(s),\bar{X}(s))|^2\\
\leq& \bar{p}^2(\sup_{0\leq r\leq t\wedge\theta_{\ell}}|e_{\D}(r)|^{\bar{p}})|e_{\D}(s)|^{\bar{p}-2}|g(s,Y(s))-g_{\D}(\bar{\rho}(s),\bar{X}(s))|^2.
\end{split}
\end{equation*}
Hence, using the inequality $(ab)^{1/2}\leq a/\l+\l b$ valid for any $a,b>0$ and $\l>0$, we have
\begin{equation}
\label{mt}
\begin{split}
&([M,M]_{t\wedge \theta_{\ell}})^{1/2}\\
\leq &\bar{p}\left(\sup_{0\leq r\leq t\wedge\theta_{\ell}}|e_{\D}(r)|^{\bar{p}}\int^{t\wedge\theta_{\ell}}_{0}|e_{\D}(s)|^{\bar{p}-2}|g(s,Y(s))-g_{\D}(\bar{\rho}(s),\bar{X}(s))|^2dE(s)\right)^{\frac{1}{2}}\\
\leq & \bar{p}\left(\frac{\sup_{0\leq r\leq t\wedge\theta_{\ell}}|e_{\D}(r)|^{\bar{p}}}{2\bar{p}b_1}+2\bar{p}b_1\int^{t\wedge\theta_{\ell}}_{0}|e_{\D}(s)|^{\bar{p}-2}|g(s,Y(s))-g_{\D}(\bar{\rho}(s),\bar{X}(s))|^2dE(s)\right).
\end{split}
\end{equation}
Combing \eqref{error}, \eqref{mt} and then apply the BDG inequality, we have
\begin{equation}\label{err2}
\begin{split}
&\EE_B\left(\sup_{0\leq t\leq T}|e_{\D}(t\wedge \theta_{\ell})|^{\bar{p}}\right)\\
\leq&b\EE_B\bigg(\sup_{0\leq t\leq T}\int^{t\wedge\theta_{\ell}}_{0}\bar{p}|e_{\D}(s)|^{\bar{p}-2}\bigg(|e_{\D}(s)|^{\mathrm{T}}\left(f(s,Y(s))-f_{\D}(\bar{\rho}(s),\bar{X}(s))\right)\\
\quad& +\frac{5\bar{p}-1}{2}|g(s,Y(s))-g_{\D}(\bar{\rho}(s),\bar{X}(s))|^2 \bigg)dE(s)\bigg)+\frac{1}{2}\sup_{0\leq r\leq t\wedge \theta_{\ell}}|e_{\D}(r)|^{\bar{p}},
\end{split}
\end{equation}
where $b=b_1 \vee 1$. Noting
\begin{equation*}
\begin{split}
&\frac{5\bar{p}-1}{2}|g(s,Y(s))-g_{\D}(\bar{\rho}(s),\bar{X}(s))|^2\\
=& \frac{5\bar{p}-1}{2}\left|g(s,Y(s))-g(s,X(s))+g(s,X(s))-g_{\D}(\bar{\rho}(s),\bar{X}(s))\right|^2\\
\leq& \frac{5\bar{p}-1}{2}\bigg(|g(s,Y(s))-g(s,X(s))|^2+2|g(s,Y(s))-g(s,X(s))|\\
& \times|g(s,X(s))-g_{\D}(\bar{\rho}(s),\bar{X}(s))|+|g(s,X(s))-g_{\D}(\bar{\rho}(s),\bar{X}(s))|^2\bigg),
\end{split}
\end{equation*}
using the Young inequality,
\begin{equation*}
\begin{split}
&2|g(s,Y(s))-g(s,X(s))||g(s,X(s))-g_{\D}(\bar{\rho}(s),\bar{X}(s))|\\
\leq& 2\bigg(\frac{1}{2}\cdot\frac{5p-5\bar{p}}{5\bar{p}-1}|g(s,Y(s))-g(s,X(s))|^2+\frac{1}{2}\cdot\frac{5\bar{p}-1}{5p-5\bar{p}}|g(s,X(s))-g_{\D}(\bar{\rho}(s),\bar{X}(s))|^2\bigg)\\
=& \frac{5p-5\bar{p}}{5\bar{p}-1}|g(s,Y(s))-g(s,X(s))|^2+\frac{5\bar{p}-1}{5p-5\bar{p}}|g(s,X(s))-g_{\D}(\bar{\rho}(s),\bar{X}(s))|^2,
\end{split}
\end{equation*}
then we can show that
\begin{equation*}
\begin{split}
&\frac{5\bar{p}-1}{2}|g(s,Y(s))-g_{\D}(\bar{\rho}(s),\bar{X}(s))|^2\\
\leq& \frac{5p-1}{2}|g(s,Y(s))-g(s,X(s))|^2+\frac{(5\bar{p}-1)(5p-1)}{2(5p-5\bar{p})}|g(s,X(s))-g_{\D}(\bar{\rho}(s),\bar{X}(s))|^2.
\end{split}
\end{equation*}
We get from \eqref{err2} that
\begin{equation}
\label{sup}
\begin{split}
\EE_B\left(\sup_{0\leq t\leq T}|e_{\D}( t\wedge\theta_{\ell})|^{\bar{p}}\right)\leq& \frac{1}{2}\sup_{0\leq r\leq t\wedge\theta_{\ell}}|e_{\D}(r)|^{\bar{p}}\\
&+\EE_B\sup_{0\leq t\leq T}[J_1]+\EE_B\sup_{0\leq t\leq T}[J_2]+\EE_B\sup_{0\leq t\leq T}[J_3],
\end{split}
\end{equation}
where
\begin{equation*}
\begin{split}
J_1&:=b\int^{t\wedge\theta_{\ell}}_{0}\bar{p}|e_{\D}(s)|^{\bar{p}-2}\bigg(e_{\D}^{\mathrm{T}}(s)\left(f(s,Y(s))-f(s,X(s))\right)\\
&\quad+\frac{5p-1}{2}|g(s,Y(s))-g(s,X(s))|^2\bigg)dE(s),
\end{split}
\end{equation*}
\begin{equation*}
\begin{split}
J_2&:=b\int^{t\wedge\theta_{\ell}}_{0}\bar{p}|e_{\D}(s)|^{\bar{p}-2}\bigg(e_{\D}^{\mathrm{T}}(s)\left(f(s,X(s))-f(\bar{\r}(s),X(s))\right)\\
&\quad+\frac{(5\bar{p}-1)(5p-1)}{(5p-5\bar{p})}|g(s,X(s))-g(\bar{\r}(s),X(s))|^2\bigg)dE(s),
\end{split}
\end{equation*}
\begin{equation*}
\begin{split}
J_3&:=b\int^{t\wedge\theta_{\ell}}_{0}\bar{p}|e_{\D}(s)|^{\bar{p}-2}\bigg(e_{\D}^{\mathrm{T}}(s)\left(f(\bar{\r}(s),X(s))-f_{\D}(\bar{\r}(s),\bar{X}(s))\right)\\
&\quad+\frac{(5\bar{p}-1)(5p-1)}{(5p-5\bar{p})}|g(\bar{\r}(s),X(s))-g_{\D}(\bar{\r}(s),\bar{X}(s))|^2\bigg)dE(s).
\end{split}
\end{equation*}
By Assumption \ref{assumption1-2}, we have
\begin{equation}
\label{J1}
J_1\leq C_1\int^{t\wedge\theta_{\ell}}_{0}|e_{\D}(s)|^{\bar{p}}dE(s),
\end{equation}
where $C_1=b\bar{p}K$.
Using the Young inequality and Assumption \ref{assumption2-4}, we can derive
\begin{equation*}
\begin{split}
J_2&\leq b\int^{t\wedge\theta_{\ell}}_{0}\bar{p}|e_{\D}(s)|^{\bar{p}-2}\bigg(\frac{1}{2}|e_{\D}(s)|^2+\frac{1}{2}|f(s,X(s))-f(\bar{\r}(s),X(s))|^2\\
&\quad+\frac{(5\bar{p}-1)(5p-1)}{(5p-5\bar{p})}|g(s,X(s))-g(\bar{\r}(s),X(s))|^2\bigg)dE(s)\\
&\leq \frac{(\bar{p}-1)(5p-5\bar{p})+(\bar{p}-2)(5\bar{p}-1)(5p-1)}{(5p-5\bar{p})}b\int^{t\wedge\theta_{\ell}}_{0}|e_{\D}(s)|^{\bar{p}}dE(s)\\
&\quad +b\int^{t\wedge\theta_{\ell}}_{0}|f(s,X(s))-f(\bar{\r}(s),X(s))|^{\bar{p}}dE(s)\\
&\quad +\frac{2(5\bar{p}-1)(5p-1)}{(5p-5\bar{p})}b\int^{t\wedge\theta_{\ell}}_{0}|g(s,X(s))-g(\bar{\r}(s),X(s))|^{\bar{p}}dE(s)\\
&\leq C_2'\bigg(\int^{t\wedge\theta_{\ell}}_{0}|e_{\D}(s)|^{\bar{p}}dE(s)+\int^{t\wedge\theta_{\ell}}_{0}H_1^{\bar{p}}(1+|X(s)|^{(\a+1)\bar{p}})\D^{\gamma_f\bar{p}}dE(s)\\
&\quad +\int^{t\wedge\theta_{\ell}}_{0}H_2^{\bar{p}}(1+|X(s)|^{(\a+1)\bar{p}})\D^{\gamma_g\bar{p}}dE(s)\bigg),
\end{split}
\end{equation*}
where
\begin{equation*}
C_2'=\max{\left\{\frac{(\bar{p}-1)(5p-5\bar{p})+(\bar{p}-2)(5\bar{p}-1)(5p-1)}{(5p-5\bar{p})},1,\frac{2(5\bar{p}-1)(5p-1)}{(5p-5\bar{p})}\right\}}b.
\end{equation*}
Then by Lemma \ref{lemma3}, we obtain
\begin{equation}
\label{J2}
J_2\leq C_2\bigg(\int^{t\wedge\theta_{\ell}}_{0}|e_{\D}(s)|^{\bar{p}}dE(s)+\D^{\gamma_f\bar{p}}E(T)+\D^{\gamma_g\bar{p}}E(T)\bigg),
\end{equation}
where
\begin{equation*}
C_2=C_2'\max\{1,H_1^{\bar{p}}(1+C)^{\a+1},H_2^{\bar{p}}(1+C)^{\a+1}\}.
\end{equation*}
Rearranging $J_3$ gives
\begin{equation}
\label{J3}
\begin{split}
J_3&\leq b\int^{t\wedge\theta_{\ell}}_{0}\bar{p}|e_{\D}(s)|^{\bar{p}-2}\bigg(e_{\D}^{\mathrm{T}}(s)\left(f(\bar{\r}(s),X(s))-f(\bar{\r}(s),\bar{X}(s))\right)\\
&\quad+\frac{2(5\bar{p}-1)(5p-1)}{(5p-5\bar{p})}|g(\bar{\r}(s),X(s))-g(\bar{\r}(s),\bar{X}(s))|^2\bigg)dE(s)\\
&\quad +b\int^{t\wedge\theta_{\ell}}_{0}\bar{p}|e_{\D}(s)|^{\bar{p}-2}\bigg(e_{\D}^{\mathrm{T}}(s)\left(f(\bar{\r}(s),\bar{X}(s))-f_{\D}(\bar{\r}(s),\bar{X}(s))\right)\\
&\quad+\frac{2(5\bar{p}-1)(5p-1)}{(5p-5\bar{p})}|g(\bar{\r}(s),\bar{X}(s))-g_{\D}(\bar{\r}(s),\bar{X}(s))|^2\bigg)dE(s)\\
&=:J_{31}+J_{32}.
\end{split}
\end{equation}
First estimate $J_{31}$. By the Young inequality and Assumption \ref{asspumption1} we can show that
\begin{equation}
\label{J311}
\begin{split}
J_{31}&\leq b\int^{t\wedge\theta_{\ell}}_{0}\bar{p}|e_{\D}(s)|^{\bar{p}-2}\bigg(\frac{1}{2}|e_{\D}(s)|^2+\frac{1}{2}|f(\bar{\r}(s),X(s))-f(\bar{\r}(s),\bar{X}(s))|^2\\
&\quad+\frac{2(5\bar{p}-1)(5p-1)}{(5p-5\bar{p})}|g(\bar{\r}(s),X(s))-g(\bar{\r}(s),\bar{X}(s))|^2\bigg)dE(s)\\
&\leq \frac{(\bar{p}-1)(5p-5\bar{p})+2(\bar{p}-2)(5\bar{p}-1)(5p-1)}{(5p-5\bar{p})}b\int^{t\wedge\theta_{\ell}}_{0}|e_{\D}(s)|^{\bar{p}}dE(s)\\
&\quad +b\int^{t\wedge\theta_{\ell}}_{0}|f(\bar{\r}(s),X(s))-f(\bar{\r}(s),\bar{X}(s))|^{\bar{p}}dE(s)\\
&\quad +\frac{4(5\bar{p}-1)(5p-1)}{(5p-5\bar{p})}b\int^{t\wedge\theta_{\ell}}_{0}|g(\bar{\r}(s),X(s))-g(\bar{\r}(s),\bar{X}(s))|^{\bar{p}}dE(s)\\
&\leq C_{31}\bigg(\int^{t\wedge\theta_{\ell}}_{0}|e_{\D}(s)|^{\bar{p}}dE(s)+2L\int^{t\wedge\theta_{\ell}}_{0}(1+|X(s)|^{\a\bar{p}}+|\bar{X}(s)|^{\a\bar{p}})|X(s)-\bar{X}(s)|^{\bar{p}}dE(s)\bigg),
\end{split}
\end{equation}
where
\begin{equation*}
C_{31}=\max\left\{\frac{(\bar{p}-1)(5p-5\bar{p})+2(\bar{p}-2)(5\bar{p}-1)(5p-1)}{(5p-5\bar{p})},1,\frac{4(5\bar{p}-1)(5p-1)}{(5p-5\bar{p})}\right\}b.
\end{equation*}
Similarly, we can show
\begin{equation*}
\begin{split}
J_{32}&\leq C_{32}\bigg(\int^{t\wedge\theta_{\ell}}_{0}|e_{\D}(s)|^{\bar{p}}dE(s)+ \int^{t\wedge\theta_{\ell}}_{0}|f(\bar{\r}(s),\bar{X}(s))-f_{\D}(\bar{\r}(s),\bar{X}(s))|^{\bar{p}}\\
&\quad +|g(\bar{\r}(s),\bar{X}(s))-g_{\D}(\bar{\r}(s),\bar{X}(s))|^{\bar{p}}dE(s)\bigg).
\end{split}
\end{equation*}
Recalling the definition of truncated EM method and Assumption \ref{asspumption1} gives
\begin{equation}
\label{J322}
\begin{split}
J_{32} &\leq C_{32}\bigg(\int^{t\wedge\theta_{\ell}}_{0}|e_{\D}(s)|^{\bar{p}}dE(s) +\int^{t\wedge\theta_{\ell}}_{0}|f(\bar{\r}(s),\bar{X}(s))-f(\bar{\r}(s),\pi_{\D}(\bar{X}(s)))|^{\bar{p}}\\
&\quad +|g(\bar{\r}(s),\bar{X}(s))-g(\bar{\r}(s),\pi_{\D}(\bar{X}(s)))|^{\bar{p}}dE(s)\bigg)\\
&\leq C_{32}\bigg(\int^{t\wedge\theta_{\ell}}_{0}|e_{\D}(s)|^{\bar{p}}dE(s)\\
&\quad +2L\int^{t\wedge\theta_{\ell}}_{0}(1+|\bar{X}(s)|^{\a\bar{p}}+|\pi_{\D}(\bar{X}(s))|^{\a\bar{p}})|\bar{X}(s)-\pi_{\D}(\bar{X}(s))|^{\bar{p}}dE(s)\bigg).
\end{split}
\end{equation}
Putting \eqref{J1}, \eqref{J2}, \eqref{J3}, \eqref{J311} and \eqref{J322} into \eqref{sup}, we get
\begin{equation}\label{edelta}
\begin{split}
\EE_B\left(\sup_{0\leq t\leq T}|e_{\D}(t\wedge\theta_{\ell})|^{\bar{p}}\right)&\leq 2C_1\EE_B\int^{T}_{0}|e_{\D}(s)|^{\bar{p}}dE(s)+2C_2\bigg(\EE_B\int^{T}_{0}|e_{\D}(s)|^{\bar{p}}dE(s)\\
&\quad+\D^{\gamma_f\bar{p}}E(T)+\D^{\gamma_g\bar{p}}E(T)\bigg)+\EE_B[I_{31}]+\EE_B[I_{32}],
\end{split}
\end{equation}
where
\begin{equation*}
I_{31}:=2C_{31}\bigg(\int^{T}_{0}|e_{\D}(s)|^{\bar{p}}dE(s)+2L\int^{T}_{0}(1+|X(s)|^{\a\bar{p}}+|\bar{X}(s)|^{\a\bar{p}})|X(s)-\bar{X}(s)|^{\bar{p}}dE(s)\bigg),
\end{equation*}
and
\begin{equation*}
\begin{split}
I_{32}:=&2C_{32}\bigg(\int^{T}_{0}|e_{\D}(s)|^{\bar{p}}dE(s)+2L\int^{T}_{0}(1+|\bar{X}(s)|^{\a\bar{p}}+|\pi_{\D}(\bar{X}(s))|^{\a\bar{p}})\\
& \times |\bar{X}(s)-\pi_{\D}(\bar{X}(s))|^{\bar{p}}dE(s)\bigg).
\end{split}
\end{equation*}
Using the H\"older inequality, Lemma \ref{lemma2} and Lemma \ref{lemma3}
\begin{equation}
\label{J31}
\begin{split}
\EE_B[I_{31}]&\leq 2C_{31}\bigg(\EE_B\int^{T}_{0}|e_{\D}(s)|^{\bar{p}}dE(s)+2L(C+1)^{\frac{\a\bar{p}}{q}}\int^{T}_{0}\left(\EE_B|X(s)-\bar{X}(s)|^q\right)^{\bar{p}/q}dE(s)\bigg)\\
&\leq 2C_{31}\bigg(\EE_B\int^{T}_{0}|e_{\D}(s)|^{\bar{p}}dE(s)+2L(C+1)^{\frac{\a\bar{p}}{q}}2^{(q-1)\bar{p}/q}\D^{\bar{p}/2}(\k(\D))^{\bar{p}}E(T)\bigg).
\end{split}
\end{equation}
Then, by the H\"older inequality and Lemma \ref{lemma3}, we obtain
\begin{equation}
\label{J32}
\begin{split}
J_{32}&\leq 2C_{32}\bigg(\EE_B\int^{T}_{0}|e_{\D}(s)|^{\bar{p}}dE(s)+2L\int^{T}_{0}\left(\EE_B(1+|\bar{X}(s)|^q+|\pi_{\D}(\bar{X}(s))|^q)\right)^{\frac{\a\bar{p}}{q}}\\
&\quad \times [\EE_B|\bar{X}(s)-\pi_{\D}(\bar{X}(s))|^\frac{q\bar{p}}{q-\a\bar{p}}]^{\frac{q-\a\bar{p}}{q}}dE(s)\bigg)\\
&\leq 2C_{32}\bigg(\EE_B\int^{T}_{0}|e_{\D}(s)|^{\bar{p}}dE(s)+4L(C+1)^{\frac{\a\bar{p}}{q}}\\
&\quad \times \int^{T}_{0}(\EE_B[I_{\{|\bar{X}(s)|>\m^{-1}(\k(\D))\}}|\bar{X}(s)|^{\frac{q\bar{p}}{q-\a\bar{p}}}])^{\frac{q-\a\bar{p}}{q}}dE(s)\bigg)\\
&\leq 2C_{32}\bigg(\EE_B\int^{T}_{0}|e_{\D}(s)|^{\bar{p}}dE(s)+4L(C+1)^{\frac{\a\bar{p}}{q}}\\
&\quad \times \int^{T}_{0}([\PP\{|\bar{X}(s)|>\m^{-1}(\k(\D))\}]^{\frac{q-(\a+1)\bar{p}}{q-\a\bar{p}}}\left(\EE_B|\bar{X}(s)|^q\right)^{\frac{\bar{p}}{q-\a\bar{p}}})^{\frac{q-\a\bar{p}}{q}}dE(s)\bigg)\\
&\leq 2C_{32}\bigg(\EE_B\int^{T}_{0}|e_{\D}(s)|^{\bar{p}}dE(s)+4L(C+1)^{\frac{(\a+1)\bar{p}}{q}}\int^{T}_{0}\left(\frac{\EE_B|\bar{X}(s)|^q}{(\m^{-1}(\k(\D)))^q}\right)^{\frac{q-(\a+1)\bar{p}}{q}}dE(s)\bigg)\\
&\leq 2C_{32}\bigg(\EE_B\int^{T}_{0}|e_{\D}(s)|^{\bar{p}}dE(s)+4L(C+1)(\m^{-1}(\k(\D)))^{(\a+1)\bar{p}-q}E(T)\bigg).
\end{split}
\end{equation}
Substituting \eqref{J31} and \eqref{J32} into \eqref{edelta} gives
\begin{equation*}
\begin{split}
\EE_B \left( \sup_{0\leq t\leq T}|e_{\D}(t)|^{\bar{p}} \right)&\leq C\bigg(\int^{T}_{0}\EE_B \left( \sup_{0\leq s\leq T}|e_{\D}(s)|^{\bar{p}}\right)dE(s)+\D^{\gamma_f\bar{p}}+\D^{\gamma_g\bar{p}}\\
&\quad +\D^{\bar{p}/2}(\k(\D))^{\bar{p}}+(\m^{-1}(\k(\D)))^{(\a+1)\bar{p}-q}\bigg).
\end{split}
\end{equation*}
An application of the Gronwall inequality yields that
\begin{equation*}
\EE_B \left( \sup_{0\leq t\leq T}|e_{\D}(t)|^{\bar{p}} \right)\leq C\bigg(\D^{\gamma_f\bar{p}}+\D^{\gamma_g\bar{p}}+\D^{\bar{p}/2}(\k(\D))^{\bar{p}}+(\m^{-1}(\k(\D)))^{(\a+1)\bar{p}-q}\bigg),
\end{equation*}
where $C$ depends on $e^{\l E(T)}$ for some $\l > 0$, but is independent from $\Delta$.
\par
Taking $\EE_D$ on both sides gives the required assertion \eqref{th321}. The other assertion \eqref{th322} follows from \eqref{th321} and Lemma \ref{lemma2}. Therefore, the proof is completed. \eproof

\section{Numerical simulations}
Two examples are considered in this section.
\begin{expl}\label{ex}
Consider a time-changed SDE
\begin{equation}\label{ex1}
\left\{
\begin{array}{lr}
dY(t)=\left([t(1-t)]^{\frac{1}{2}}Y^2(t)-2Y^5(t)\right)dE(t)+\left([t(1-t)]^{\frac{1}{4}}Y^2(t)\right)dB(E(t)),&\\
Y(0)=2,&
\end{array}
\right.
\end{equation}
with  $T=1$.
\end{expl}
\par
For any $p>1$, we can see
\begin{equation*}
\begin{split}
&(x-y)^{\mathrm{T}}(f(t,x)-f(t,y))+\frac{5p-1}{2}|g(t,x)-g(t,y)|^2\\
\leq& (x-y)^2\bigg([t(1-t)]^{\frac{1}{2}}(x+y)-x^4-y^4+(5p-1)[t(1-t)]^{\frac{1}{2}}(x^2+y^2)\bigg)\\
\leq& C(x-y)^2,
\end{split}
\end{equation*}
where the Young inequality is used. Note that the last inequality is due to the fact that polynomials with the negative coefficients for the highest order term can always be bounded from above. This indicates that Assumption \ref{assumption1-2} holds.
\par
In the similar manner, for any $q>1$ and any $t\in [0,1]$, we have
\begin{equation*}
\begin{split}
& x^{\mathrm{T}}f(t,x)+\frac{5q-1}{2}|g(t,x)|^2\\
=& [t(1-t)]^{\frac{1}{2}}x^3-x^6+\frac{5q-1}{2}[t(1-t)]^{\frac{1}{2}}x^4\\
\leq & C(1+|x|^2),
\end{split}
\end{equation*}
which means that Assumption \ref{assumption2} is satisfied.
\par
Using the mean theorem for the temporal variable, Assumption \ref{asspumption1} and Assumption \ref{assumption2-4} are satisfied with $\a=4$, $\gamma_f=\frac{1}{2}$, $\gamma_f=\frac{1}{4}$. According to Theorem \ref{theorem3-2}, we know that
\begin{equation*}
\EE \left( \sup_{0\leq t\leq T}|Y(t)-X(t)|^{\bar{p}} \right) \leq C\bigg(\D^{\frac{\bar{p}}{2}}+\D^{\frac{\bar{p}}{4}}+\D^{\frac{\bar{p}}{2}}(\k(\D))^{\bar{p}}+(\mu^{-1}(\k(\D)))^{5\bar{p}-q}\bigg),
\end{equation*}
and
\begin{equation*}
\EE \left( \sup_{0\leq t\leq T}|Y(t)-\bar{X}(t)|^{\bar{p}} \right) \leq C\bigg(\D^{\frac{\bar{p}}{2}}+\D^{\frac{\bar{p}}{4}}+\D^{\frac{\bar{p}}{2}}(\k(\D))^{\bar{p}}+(\mu^{-1}(\k(\D)))^{5\bar{p}-q}\bigg).
\end{equation*}
Due to that
\begin{equation*}
 \sup_{0\leq t\leq T}\sup_{|x|\leq u}(|f(t,x)|\vee |g(t,x)|)\leq 3u^5,\quad  \forall u\geq 1,
\end{equation*}
we choose $\m (u)=3u^5$ and $\k(\D)=\D^{-\varepsilon}$, for any $\varepsilon \in (0,1/4]$. As a result, $\m^{-1}(u)=\left(u/3\right)^{1/5}$ and $\m^{-1}(\k(\D))=\left(\D^{-\varepsilon}/3\right)^{1/5}$. Choosing $p$ sufficiently large, we can get from Theorem \ref{theorem3-1} that
\begin{equation*}
\EE \left( \sup_{0\leq t\leq 1}|Y(t)-X(t)|^{\bar{p}} \right) \leq C\D^{\bar{p}/4},
\end{equation*}
and
\begin{equation*}
\EE \left( \sup_{0\leq t\leq 1}|Y(t)-\bar{X}(t)|^{\bar{p}} \right) \leq C\D^{\bar{p}/4},
\end{equation*}
which imply that the convergence rate of truncated EM method for the time-change SDE \eqref{ex1} is $1/4$.
\par
Let us compute the approximation of the mean square error. We run M=100 independent trajectories using \eqref{eq:DiscreteEM} for every different step sizes, $10^{-2}$, $10^{-3}$, $10^{-4}$, $10^{-6}$. Because it is hard to find the true solution for the SDE, the numerical solution with the step size $10^{-6}$ is regarded as the exact solution.
\begin{figure}[htph]
\centering
\includegraphics[width=0.80\textwidth]{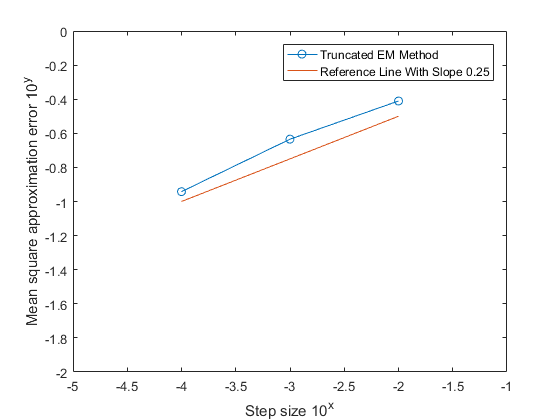}
\caption{Convergence rate of Example \ref{ex}}
\end{figure}
It is not hard to see from Figure 1 that the strong convergence rate is approximately 0.25. To see it more clearly, applying the linear regression shows that the slope of the line of errors is about 0.26581, which is in line with the theoretical result.
\begin{expl}\label{exp}
Consider a two-dimensional time-changed SDE
\begin{equation}\label{ex2}
\left\{
\begin{array}{lr}
dx_1(t)=\left([(t-1)(2-t)]^{\frac{1}{5}}x_1^2(t)-2x_2^5(t)\right)dE(t)+\left([(t-1)(2-t)]^{\frac{2}{5}}x_2^2(t)\right)dB(E(t)),&\\
dx_2(t)=\left([(t-1)(2-t)]^{\frac{1}{5}}x_2^2(t)-2x_1^5(t)\right)dE(t)+\left([(t-1)(2-t)]^{\frac{2}{5}}x_1^2(t)\right)dB(E(t)).&
\end{array}
\right.
\end{equation}
\end{expl}
It is clear that
\begin{equation*}
f(t,x)=
\begin{pmatrix}
[(t-1)(2-t)]^{\frac{1}{5}}x_1^2-2x_2^5\\
[(t-1)(2-t)]^{\frac{1}{5}}x_2^2-2x_1^5
\end{pmatrix}
~~~~~\text{and}~~~~~
g(t,x)=
\begin{pmatrix}
[(t-1)(2-t)]^{\frac{2}{5}}x_2^2\\
[(t-1)(2-t)]^{\frac{2}{5}}x_1^2
\end{pmatrix}
.
\end{equation*}
\par
For any $x,y\in\RR$, it is easy to show that
\begin{equation*}
\begin{split}
&(x-y)^{\mathrm{T}}(f(t,x)-f(t,y))+\frac{5p-1}{2}|g(t,x)-g(t,y)|^2\\
\leq & (x_1-y_1)^2\left([(t-1)(2-t)]^{\frac{1}{5}}(x_1+y_1)-2(x_2^4+x_2^3y_2+x_2^2y_2^2+x_2y_2^3+y_2^4)\right)\\
\quad &+(x_2-y_2)^2\left([(t-1)(2-t)]^{\frac{1}{5}}(x_2+y_2)-2(x_1^4+x_1^3y_1+x_1^2y_1^2+x_1y_1^3+y_1^4)\right)\\
\quad &+\frac{5p-1}{2}\bigg([(t-1)(2-t)]^{\frac{4}{5}}(x_2^2-y_2^2)^2+2[(t-1)(2-t)]^{\frac{2}{5}}(x_2^2-y_2^2)(x_1^2-y_1^2)\\
\quad &+[(t-1)(2-t)]^{\frac{4}{5}}(x_1^2-y_1^2)^2\bigg)\\
\leq & (x_1-y_1)^2\left([(t-1)(2-t)]^{\frac{1}{5}}(x_1+y_1)-2(x_2^4+x_2^3y_2+x_2^2y_2^2+x_2y_2^3+y_2^4)\right)\\
\quad &+(x_2-y_2)^2\left([(t-1)(2-t)]^{\frac{1}{5}}(x_2+y_2)-2(x_1^4+x_1^3y_1+x_1^2y_1^2+x_1y_1^3+y_1^4)\right)\\
\quad &+(x_2-y_2)^2\frac{(5p-1)}{2}\bigg([(t-1)(2-t)]^{\frac{4}{5}}(x_2+y_2)^2+(x_1-y_1)(x_2-y_2)(5p-1)\\
&\quad \times [(t-1)(2-t)]^{\frac{2}{5}}+(x_1-y_1)^2\frac{(5p-1)}{2}[(t-1)(2-t)]^{\frac{4}{5}}(x_1+y_1)^2\bigg).
\end{split}
\end{equation*}
But
\begin{equation*}
-2(x^3y+xy^3)=-2xy(x^2+y^2)\leq (x^2+y^2)^2=x^4+y^4+2x^2y^2.
\end{equation*}
Therefore, for any $t\in[0,1]$
\begin{equation*}
\begin{split}
&(x-y)^{\mathrm{T}}(f(t,x)-f(t,y))+\frac{5p-1}{2}|g(t,x)-g(t,y)|^2\\
\leq &(x_1-y_1)^2\bigg([(t-1)(2-t)]^{\frac{1}{5}}(x_1+y_1)-x_2^4-y_2^4+(5p-1)[(t-1)(2-t)]^{\frac{4}{5}}(x_1^2+y_1^2)\\
\quad &+\frac{5p-1}{2}[(t-1)(2-t)]^{\frac{2}{5}}\bigg)+(x_2-y_2)^2\bigg([(t-1)(2-t)]^{\frac{1}{5}}(x_2+y_2)-x_1^4-y_1^4\\
\quad &+(5p-1)[(t-1)(2-t)]^{\frac{4}{5}}(x_2^2+y_2^2)+\frac{5p-1}{2}[(t-1)(2-t)]^{\frac{2}{5}}\bigg)\\
\leq &C(x-y)^2,
\end{split}
\end{equation*}
where the last inequality is due to the fact that polynomials with the negative coefficients for the highest order term can always be bounds from above. This indicates that Assumption \ref{assumption1-2} holds.
\par
In the similar manner, for any $q>1$ and any $t\in[0,1]$, we have
\begin{equation*}
\begin{split}
& x^{\mathrm{T}}f(t,x)+\frac{5q-1}{2}|g(t,x)|^2\\
=& ([(t-1)(2-t)]^{\frac{1}{5}}x_1^3-2x_1x_2^5)+([(t-1)(2-t)]^{\frac{1}{5}}x_2^3-2x_2x_1^5)\\
\quad &+\frac{5q-1}{2}|[(t-1)(2-t)]^{\frac{2}{5}}x_1x_2(x_2+x_1)|^2\\
\leq & [(t-1)(2-t)]^{\frac{1}{5}}(x_1^3+x_2^3)-2x_1x_2(x_1^4+x_2^4)+(5q-1)[(t-1)(2-t)]^{\frac{4}{5}}(x_1^6+x_2^6)\\
\leq & C(1+|x|^2),
\end{split}
\end{equation*}
which means that Assumption \ref{assumption2} is satisfied. Then, we consider Assumption \ref{asspumption1}, by the means theorem and the elementary inequality we can have
\begin{equation*}
\begin{split}
&|f(t,x)-f(t,y)|\\
\leq& |[(t-1)(2-t)]^{\frac{1}{5}}(x_1^2-y_1^2)-2(x_2^5-y_2^5)+[(t-1)(2-t)]^{\frac{1}{5}}(x_2^2-y_2^2)-2(x_1^5-y_1^5)|\\
\leq& (x_1-y_1)|[(t-1)(2-t)]^{\frac{1}{5}}(x_1+y_1)-x_1^4-y_1^4|+(x_2-y_2)|[(t-1)(2-t)]^{\frac{1}{5}}\\
&\times (x_2+y_2)-x_2^4-y_2^4|\\
\leq& C(1+|x|^4+|y|^4)|x-y|.
\end{split}
\end{equation*}
Similarly, we derive that
\begin{equation*}
|g(t,x)-g(t,y)|\leq C(1+|x|+|y|)|x-y|.
\end{equation*}
Then, assume that $\gamma_f\in (0,1]$ and $\gamma_g\in (0,1]$, for any $s,t\in[0,T]$, using the mean value theorem for the temporal variable
\begin{equation*}
\begin{split}
&|f(s,x)-f(t,x)|\\
\leq& |([(s-1)(2-s)]^{\frac{1}{5}}-[(t-1)(2-t)]^{\frac{1}{5}})x_1^2+([(s-1)(2-s)]^{\frac{1}{5}}-[(t-1)(2-t)]^{\frac{1}{5}})x_2^2|\\
\leq& C_1|s-t|^{\frac{1}{5}}x_1^2+C_2|s-t|^{\frac{1}{5}}x_2^2,
\end{split}
\end{equation*}
and
\begin{equation*}
\begin{split}
&|g(s,x)-g(t,x)|\\
\leq& |([(s-1)(2-s)]^{\frac{2}{5}}-[(t-1)(2-t)]^{\frac{2}{5}})x_2^2+([(s-1)(2-s)]^{\frac{2}{5}}-[(t-1)(2-t)]^{\frac{2}{5}})x_1^2|\\
\leq& C_1|s-t|^{\frac{2}{5}}x_2^2+C_2|s-t|^{\frac{2}{5}}x_1^2.
\end{split}
\end{equation*}
Thus, Assumptions \ref{asspumption1} and \ref{assumption2-4} are satisfied with $\a=4$, $\gamma_f=1/5$ and $\gamma_g=2/5$. According to Theorem \ref{theorem3-2}, we know that
\begin{equation*}
\EE \left( \sup_{0\leq t\leq T}|Y(t)-X(t)|^{\bar{p}} \right) \leq C\bigg(\D^{\frac{\bar{p}}{5}}+\D^{\frac{2\bar{p}}{5}}+\D^{\frac{\bar{p}}{2}}(\k(\D))^{\bar{p}}+(\mu^{-1}(\k(\D)))^{5\bar{p}-q}\bigg),
\end{equation*}
and
\begin{equation*}
\EE \left( \sup_{0\leq t\leq T}|Y(t)-\bar{X}(t)|^{\bar{p}}\right)\leq C\bigg(\D^{\frac{\bar{p}}{5}}+\D^{\frac{2\bar{p}}{5}}+\D^{\frac{\bar{p}}{2}}(\k(\D))^{\bar{p}}+(\mu^{-1}(\k(\D)))^{5\bar{p}-q}\bigg).
\end{equation*}
The remainder of the proof is omitted since it is same to the argument given in the last part of the Example \ref{ex}. Finally, we can get that
\begin{equation*}
\EE \left( \sup_{0\leq t\leq 1}|Y(t)-X(t)|^{\bar{p}} \right) \leq C\D^{\bar{p}/5},
\end{equation*}
and
\begin{equation*}
\EE \left( \sup_{0\leq t\leq 1}|Y(t)-\bar{X}(t)|^{\bar{p}} \right) \leq C\D^{\bar{p}/5},
\end{equation*}
which imply that the convergence rate of truncated EM method for the time-change SDE \eqref{ex2} is $1/5$.
\par
We run M=100 independent trajectories for every different step sizes, $10^{-2}$, $10^{-3}$, $10^{-4}$, $10^{-6}$. Because it is hard to find the true solution for the SDE, the numerical solution with the step size $10^{-6}$ is regarded as the exact solution.
\begin{figure}[htph]
\centering
\includegraphics[width=0.80\textwidth]{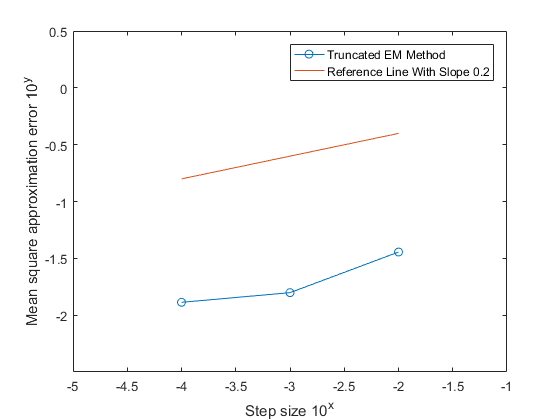}
\caption{Convergence rate of Example \ref{exp}}
\end{figure}
It is not hard to see from Figure 2 that the strong convergence rate is approximately 0.2. To see it more clearly, applying the linear regression shows that the slope of the line of errors is about 0.22123, which is in line with the theoretical result.

\section{Conclusion and future research}
In this paper, the truncated EM method is proposed for a class of time-changed SDEs and the strong convergence rate in the finite time interval is obtained. The key difference of this paper and those existing works is that both the coefficients are allowed to grow super-linearly in terms of the state variable.
\par
The asymptotic behaviour of the numerical methods is also a very interesting topic in the field of the numerical analysis for time-changed SDEs. We will discuss the asymptotic stability of the truncated EM method for time-changed SDEs in our future work.

\end{document}